# On the principal blocks of $S_{15}$ over a field of characteristic 5

Jinbao Yao

School of Mathematics and Statistics, Hubei University, Wuhan, 430062, China

E-mail: 202221104011315@stu.hubu.edu.cn

**Abstract.** In this paper, we construct the Ext-quiver of a principal block of $kS_{15}$, where $k$ is an algebraically closed field of characteristic 5.

**Keywords:** symmetric groups; Loewy structures; Ext-quiver.

## 1.Introduction

Blocks of symmetric group algebras in odd characteristic has been studied extensively with the aim of understanding algebras of wild representation type. In [1,17], Tan found the principal blocks of $kS_9$, $kS_{10}$ and $kS_{11}$ where $k$ is a field of characteristic three. In [19], Fayers found quivers of blocks (3,1) and $(3,1^2)$, the cores of the $kS_{13}$ and $kS_{14}$ blocks, respectively. However, both of them studied based on fields of characteristic three, and we conjectured whether there is a similar property for fields of larger characteristic? So we studied the quiver structure of the principal block of $kS_{15}$ with a defect of 3 when the characteristic of the fields is 5, and found that there are indeed similar properties:

  1. The simple modules do not self-extend.

  2. The $\text{Ext}^1$-space between two simple modules is at most one-dimensional.

  3. The Ext-quiver is bipartite.

## 2.Preliminaries

In this section, we give a brief account of the representation theory which we require. For more detailed accounts, we refer the reader to [2,3] for the representation theory of symmetric groups, and [4,5] for finite theory of group representations.

**Definition 2.1.** *[6] Let $\lambda = (\lambda_1, \lambda_2, \ldots, \lambda_k)$ be a partition of n (with $\lambda_i \geq \lambda_{i+1}$), $P=[\lambda_1, \lambda_2, \ldots, \lambda_k]$ is a $\lambda$-diagram and p be a characteristic of the field. Starting from the rightmost point of the first row of P, take the first p points of the rim of P from top to bottom, left to right, and these p points are called the first p-segment; the next p-segment similarly takes the p points of the rim of P in the same direction starting from the rightmost point of the next row of the row where the end of the previous p-segment is located; and when the leftmost point of the bottom row is taken point, all p-segments are taken, and notice that the number of points in the last p-segment that satisfy the requirement may be less than p.*

*For a $\lambda$-diagram P, the set of all its p-segments is called the p-edge of P. In general, we denote by $I(P)$ the diagram after removing the p-edge from the P and by $e(P)$ the number of points in the p-edge.*

*Notice.* To illustrate how the above definition relates to the partitions, we use the above definition in such a way that the default $\lambda$-diagram is the partition $\lambda$ itself.

*Example.* If $\lambda = (7,6,3,3,2,1)$ and $p=5$, then there are only three *p*-segments of the $\lambda$-diagram *P* and $e(P)=11$ and $I(P)=(5,3,2,1)$.

**Definition 2.2.** *[6] Given the characteristics and the partition P, $P' = I(P)$, we continue the process of removing the p-edge by first obtaining $I(P')$ and then continuing until obtaining the empty partition, assuming that $\alpha$ steps are taken, then we have $\alpha+1$ partitions $P_\alpha, \ldots, P_0$, where $P_\alpha = P$, $P_0 = \emptyset$, $P_{i-1} = I(P_i)$, $i=1,2,\ldots,\alpha$, and define $a_i = e(P_i)$. Suppose that $r_i$ is the number of rows of $P_i$, and for each $P_i$, define a sequence $(s_1, \ldots, s_\alpha)$ by the following formula:*

$$s_i = a_i - r_i + \varepsilon_i ,$$

*where $\varepsilon_i = 0$ if p divides $a_i$ and is 1 otherwise.*

**Lemma 2.3.** *[6,THEOREM 4.1] Let p be a characteristic of the field and P be a p-regular partition of n, and the finite sequences $(a_1, \ldots, a_\alpha)$ and $(s_1, \ldots, s_\alpha)$ be given by the above definition, then these two sequences determine a unique sequence of p-regular partitions $(Q_1, \ldots, Q_\alpha)$ such that each $Q_i$ has $s_i$ rows and $e(Q_i) = a_i$.*

**Lemma 2.4.** *[7,8] Let S be the set consisting of all p-regular partitions of n and construct a correspondence law $b_n: S \to S$, $P \mapsto Q$, where $Q$ is the p-regular partition $Q_\alpha$ obtained by Lemma 2.3, then $b_n$ is a bijection and $D^P \otimes sgn = D^{b_n(P)}$.*

**Definition 2.5.** *[10] Let $\lambda=(\lambda_1^{a_1}, \lambda_2^{a_2}, \ldots, \lambda_k^{a_k})$ be a partition of n (with $\lambda_i \geq \lambda_{i+1}$). For $1 \leq i \leq j \leq k$, we put*

$$\beta(i,j) = \lambda_i - \lambda_j + \sum_{t=i}^{j} a_t;$$

$$\gamma(i,j) = \lambda_i - \lambda_j + \sum_{t=i+1}^{j} a_t;$$

$$\lambda(i) = (\lambda_1^{a_1}, \lambda_2^{a_2}, \ldots, \lambda_{i-1}^{a_{i-1}}, \lambda_i^{a_i-1}, \lambda_i - 1, \lambda_{i+1}^{a_{i+1}}, \ldots, \lambda_k^{a_k}).$$

**Definition 2.6.** *[10, Definition 0.3] Let $\lambda = (\lambda_1^{a_1}, \lambda_2^{a_2}, \ldots, \lambda_k^{a_k})$ be a partition of n, for $i \in \{1,2,\ldots,k\}$, denote $M_i = \{j \mid 1 \leq j < i, \beta(j, i) \equiv 0 \pmod p\}$, and suppose that for all $j \in M_i$, there exists a corresponding integer $d(j)$ and $j < d(j) < i$ such that:*

  *(1) $\beta(j, d(j)) \equiv 0 \pmod p$;*

  *(2) $|\{d(j) \mid j \in M_i\}| = |M_i|$;*

*Then we call i normal.*

**Lemma 2.7.** *[10] The following statements hold:*

  *(1) 1 is normal;*

  *(2) 2 is normal if and only if $\beta(1, 2) \not\equiv 0 \pmod p$;*

  *(3) i is not normal if $i > 1$ and $\beta(i-1, i) \equiv 0 \pmod p$.*

*Proof.* (1) When $i=1$, $M_1 = \emptyset$ and thus 1 is normal by Definition 2.6. (2) In case $i=2$, $M_2 = \{1\}$ is true if $\beta(1,2) \equiv 0 \pmod p$. 2 is not normal since $d(1)$ satisfying $1 < d(1) < 2$ does not exist; on the other hand, if $\beta(1,2) \not\equiv 0 \pmod p$, then $M_2 = \emptyset$, and 2 is obviously normal. (3) $1 \leq i-1 < i$ for $i > 1$, $i-1 \in M_i$ if $\beta(i-1,i) \equiv 0 \pmod p$. Let $j = i-1$, then $i$ is not normal as $d(j)$ fulfilling $i-1 < d(j) < i$ does not exist.

**Lemma 2.8.** *[9,10,11] Let $\lambda=(\lambda_1^{a_1}, \lambda_2^{a_2}, \ldots, \lambda_k^{a_k})$ be a p-regular partition of n. Then*

$$D^\lambda \downarrow_{S_{n-1}} \cong \bigoplus \{D^{\lambda(i)} \mid i \text{ is normal}\}.$$

**Definition 2.9.** *Let p be a characteristic of the field k, B be a principal block in the symmetric group $S_{3p}$ with defect 3, and $B_s$ be a principal block in the symmetric group $S_{3p-1}$ with defect 2 and p-core $(s-1, 1^{p-s})$, where $2 \leq s \leq p$, and $B_1$ be a principal block in $S_{3p-1}$ with defect 1 and p-core $(p, 1^{p-1})$. Following these notations below, we consider the induction and restriction between the specht modules of these principal blocks.*

**Lemma 2.10.** *[12, LEMMA 3.2] Let $\tilde{\lambda}$ be a partition of $(3p-1)$ and the corresponding Specht module belongs to block $B_s$. Then there exist two partitions $\lambda$ and $\mu$ of 3p such that:*

$$S^{\tilde{\lambda}} \uparrow_{B_s}^{B} \sim S^{\lambda} + S^{\mu}, \; S^{\tilde{\lambda}} = S^{\lambda} \downarrow_{B_s} = S^{\mu} \downarrow_{B_s}.$$

*Proof.* The p-core $(s-1, 1^{p-s})$ of $B_s$ is denoted by the β-set of 3p numbers as $\{3p+s-2, 3p-1, \ldots, 2p+s, 2p+s-2, \ldots, 0\}$, when this β-set is displayed on an abacus with p runners, then there are 4 beads on runner $s-1$, 2 beads on runner $s$, and 3 beads on all the rest of the runners, so the partition corresponding to the Specht module belonging to $B_s$ can be displayed with $\langle 3^{s-2}, 4, 2, 3^{p-2} \rangle$ notation. If $\tilde{\lambda}$ is shown in $\langle 3^{s-2}, 4, 2, 3^{p-2} \rangle$ notation, the abacus display of the partition $\rho$ corresponding to the factor $S^{\rho}$ of $S^{\tilde{\lambda}} \uparrow^{B}$ can be obtained by moving the bead on the runner $s-1$ to the right on the runner $s$ in the abacus display of $\tilde{\lambda}$, and there are two beads that can be moved on the runner $s-1$. Thus there are two partitions $\lambda$ and $\mu$, such that:

$$S^{\tilde{\lambda}} \uparrow_{B_s}^{B} \sim S^{\lambda} + S^{\mu}, \; S^{\tilde{\lambda}} = S^{\lambda} \downarrow_{B_s} = S^{\mu} \downarrow_{B_s}.$$

**Definition 2.11.** *[12] If $\{\tilde{\lambda}_1, \tilde{\lambda}_2, \ldots, \tilde{\lambda}_a\}$ is the complete set of partitions whose Specht modules $S^{\tilde{\lambda}_i}$ belong to $B_s$, where $2 \leq s \leq p$, $\{\lambda_1, \mu_1, \lambda_2, \mu_2, \ldots, \lambda_a, \mu_a\}$ is the set of partitions whose Specht modules all belong to B, and with*

$$S^{\tilde{\lambda}_i} = S^{\lambda_i} \downarrow_{B_s} = S^{\mu_i} \downarrow_{B_s}, \; \lambda_i > \mu_i.$$

*Let*

$$\Lambda_s = \{\lambda_1, \lambda_2, \ldots, \lambda_a\}, \; M_s = \{\mu_1, \mu_2, \ldots, \mu_a\},$$

$$\Theta_s : \lambda_i \mapsto \tilde{\lambda}_i, \; \text{where } i \in \{1, 2, \ldots, a\}.$$

**Lemma 2.12.** *[12, LEMMA 3.6 and LEMMA 3.7] The map $\Theta_s$ has the following properties:*

(1) $\Theta_s$ *maintains the lexicographic ordering of the partition;*

(2) *The partition $\lambda_i \in \Lambda_s$ is a p-regular partition if and only if $\Theta_s(\lambda_i)$ is a p-regular partition.*

**Lemma 2.13.** *[12, THEOREM 3.9] Let $\lambda \in \Lambda_s, \mu \in M_s, \tilde{\lambda} = \Theta_s(\lambda)$, then*

(1) $D^\lambda \downarrow_{B_s} = D^{\tilde{\lambda}}, D^\mu \downarrow_{B_s} = 0;$

(2) $D^{\tilde{\lambda}} \uparrow^B$ *has two composition factors isomorphic to $D^\lambda$ and one composition factor isomorphic to $D^\mu$, where $S^\mu \downarrow_{B_s} = S^{\tilde{\lambda}}$ For other irreducible composition factors $D^\rho$ of $D^{\tilde{\lambda}} \uparrow^B$, we have $\rho \notin \Lambda_s$.*

*Proof.* If $\lambda_1 > \lambda_2 > \cdots > \lambda_a$ are all p-regular partitions in $\Lambda_s$, it follows from Lemma 2.12 that $\tilde{\lambda}_1 > \tilde{\lambda}_2 > \cdots > \tilde{\lambda}_a$ are the complete set of p-regular partitions corresponding to Specht modules belonging to $B_s$. We take any $\lambda \in \Lambda_s$ and let $\tilde{\lambda} = \Theta_s(\lambda)$, then we can assume

$$S^{\tilde{\lambda}} \sim \sum_{j=1}^{a} b_j^{\tilde{\lambda}} D^{\tilde{\lambda}_j}, b_j^{\tilde{\lambda}} \in N.$$

If $\tilde{\lambda} > \tilde{\lambda}_j$, then $b_j^{\tilde{\lambda}} = 0$; if $\tilde{\lambda} = \tilde{\lambda}_j$, then $b_j^{\tilde{\lambda}} = 1$. Now by Lemma 2.10 we have

$$S^{\tilde{\lambda}} \uparrow^B \downarrow_{B_s} \sim 2 S^{\tilde{\lambda}}.$$

Since $S^{\tilde{\lambda}_1} = D^{\tilde{\lambda}_1}$, it follows that

$$D^{\tilde{\lambda}_1} \uparrow^B \downarrow_{B_s} \sim 2 D^{\tilde{\lambda}_1}.$$

Below we use induction to show that for all $1 \leq j \leq a$, there is:

$$D^{\tilde{\lambda}_j} \uparrow^B \downarrow_{B_s} \sim 2 D^{\tilde{\lambda}_j}.$$

Assuming that for $j < l \leq a$, there is

$$D^{\tilde{\lambda}_j} \uparrow^B \downarrow_{B_s} \sim 2 D^{\tilde{\lambda}_j},$$

then

$$S^{\tilde{\lambda}_l} \uparrow^B \downarrow_{B_s} \sim [\sum_{j=1}^{a} b_j^{\tilde{\lambda}_l} D^{\tilde{\lambda}_j}] \uparrow^B \downarrow_{B_s}$$

$$\sim [\sum_{j=1}^{l-1} b_j^{\tilde{\lambda}_l} D^{\tilde{\lambda}_j} + D^{\tilde{\lambda}_l}] \uparrow^B \downarrow_{B_s}$$

$$\sim 2 \sum_{j=1}^{l-1} b_j^{\tilde{\lambda}_l} D^{\tilde{\lambda}_j} + D^{\tilde{\lambda}_l} \uparrow^B \downarrow_{B_s}.$$

Again, since

$$S^{\tilde{\lambda}}\uparrow^B\downarrow_{B_s} \sim 2S^{\tilde{\lambda}}, \quad S^{\tilde{\lambda_l}} \sim \sum_{j=1}^{a} b_j^{\tilde{\lambda_l}} D^{\tilde{\lambda_j}} \sim (\sum_{j=1}^{a-1} b_j^{\tilde{\lambda_l}} D^{\tilde{\lambda_j}} + D^{\tilde{\lambda_l}}),$$

it follows that

$$D^{\tilde{\lambda_l}}\uparrow^B\downarrow_{B_s} \sim 2D^{\tilde{\lambda_l}},$$

It is easy to show that $S^{\mu_j}$ has a composition factor that is isomorphic to $D^{\lambda_j}$, so that

$$D^{\lambda_j}\downarrow_{B_s} = D^{\tilde{\lambda_j}}, D^{\mu_j}\downarrow_{B_s} = 0.$$

Therefore, $D^{\tilde{\lambda_j}}\uparrow^B$ has two composition factors isomorphic to $D^{\lambda_j}$ and one factor isomorphic to $D^{\mu_j}$, and has no composition factors isomorphic to $D^{\lambda_i}$, where $i \neq j$. So (1) and (2) are proved.

**Lemma 2.14.** *[12, LEMMA 3.11] Let $B_1$ be a principal block with defect 1 and p-core $(p, 1^{p-1})$. Let $\tilde{\lambda}$ be a partition of $(3p-1)$ and the corresponding Specht module belongs to block $B_1$. Then there exist three partitions $\lambda, \mu$ and $v$ of $3p$ such that:*

$$S^{\tilde{\lambda}}\uparrow^B_{B_1} \sim S^\lambda + S^\mu + S^v, \quad S^{\tilde{\lambda}} = S^\lambda\downarrow_{B_1} = S^\mu\downarrow_{B_1} = S^v\downarrow_{B_1}.$$

*Proof.* The proof here is similar to that of Lemma 2.10. Clearly $\tilde{\lambda}$ is denoted by $\langle 2, 3^{p-2}, 4 \rangle$ notation, and thus the abacus display of the partition $\rho$ corresponding to the composition factor $D^\rho$ of $S^{\tilde{\lambda}}\uparrow^B$ can be obtained from the abacus display of $\tilde{\lambda}$ in which a bead on runner p is moved to the left to runner 1, and there are 3 beads on runner $p$ that can be moved, so we can obtain 3 partitions $\lambda$, $\mu$ and $v$ such that:

$$S^{\tilde{\lambda}}\uparrow^B_{B_1} \sim S^\lambda + S^\mu + S^v, \quad S^{\tilde{\lambda}} = S^\lambda\downarrow_{B_1} = S^\mu\downarrow_{B_1} = S^v\downarrow_{B_1}.$$

**Definition 2.15.** *[12] If $\{\tilde{\lambda}_1, \tilde{\lambda}_2, \ldots, \tilde{\lambda}_a\}$ is the complete set of partitions whose Specht modules $S^{\tilde{\lambda_i}}$ belong to $B_1$, $\{\lambda_1, \mu_1, v_1, \ldots, \lambda_a, \mu_a, v_a\}$ is the set of partitions whose Specht modules all belong to B, and with*

$$S^{\tilde{\lambda_i}} = S^{\lambda_i}\downarrow_{B_1} = S^{\mu_i}\downarrow_{B_1} = S^{v_i}\downarrow_{B_1}, \quad \lambda_i > \mu_i > v_i.$$

*Let*

$$\Lambda_1 = \{\lambda_1,\ldots,\lambda_a\},\ M_1 = \{\mu_1, \upsilon_1,\ldots,\mu_a, \upsilon_a\},$$

$$\Theta_1 : \lambda_i \mapsto \tilde{\lambda}_i,\ \text{where } i \in \{1,2,\ldots,a\}.$$

**Lemma 2.16.** *[12, LEMMA 3.15] The map $\Theta_1$ has the following properties:*

(1) $\Theta_1$ *maintains the lexicographic ordering of the partition;*

(2) *The partition $\lambda_i \in \Lambda_1$ is a p-regular partition if and only if $\Theta_1(\lambda_i)$ is a p-regular partition.*

**Lemma 2.17.** *[12, THEOREM 3.17] Let $\lambda \in \Lambda_1$, $\mu$, $\upsilon \in M_1$, $\tilde{\lambda} = \Theta_1(\lambda)$, then*

(1) $D^\lambda\!\downarrow_{B_1} = D^{\tilde{\lambda}}, D^\mu\!\downarrow_{B_1} = D^\upsilon\!\downarrow_{B_1} = 0;$

(2) $D^{\tilde{\lambda}}\!\uparrow^B$ *has three composition factors isomorphic to $D^\lambda$, one composition factor isomorphic to $D^\mu$, and one composition factor isomorphic to $D^\upsilon$, where*

$$S^\mu\!\downarrow_{B_1} = S^\upsilon\!\downarrow_{B_1} = S^{\tilde{\lambda}}.$$

*For other irreducible composition factors $D^\rho$ of $D^{\tilde{\lambda}}\!\uparrow^B$, we have $\rho \notin \Lambda_1$.*

**Lemma 2.18.** *[1, Lemma 2.5] Let B be a block of $kS_n$, where k is a field of any characteristic, and $\tilde{B}$ is a block of $kS_{n-1}$. Assuming $Ext^1_B(D^\lambda, D^\mu) \neq 0$ and $D^{\tilde{\lambda}} \cong soc(D^\lambda\!\downarrow_{\tilde{B}})$, one of the following statements must hold:*

(1) $Ext^1(D^{\tilde{\lambda}}\!\uparrow^B, D^\mu) = 0$; *in this case, $D^\mu$ lies in the second layer of the Loewy structure of $D^{\tilde{\lambda}}\!\uparrow^B$ and the dimension of $Ext^1(D^\lambda, D^\mu)$ is the number of copies $D^\mu$ appears in the second layer of the Loewy structure of $D^{\tilde{\lambda}}\!\uparrow^B$.*

(2) $Ext^1(D^{\tilde{\lambda}}, D^\mu\!\downarrow_{\tilde{B}}) \neq 0$; *in particular, $D^\mu\!\downarrow_{\tilde{B}} \neq 0$.*

## 3.The Ext-quiver

In this section, we let the characteristic $p=5$, so by Definition 2.9, $B$ is the principal block of $S_{15}$ and $B_s(2 \leq s \leq p)$ and $B_1$ are the principal blocks of $S_{14}$.

**Definition 3.1.** *let $\lambda$ be a partition of 15 and the weight of $\lambda$ be 3, then $\lambda$ can be represented on an abacus with 15 beads, and we make the following notation: if the abacus display of $\lambda$ with 15 beads*

*has*

 (1) *one bead of weight 3 on column i, then denote λ by* ⟨i⟩;

 (2) *one bead of weight 2 on column i and one bead of weight 1 on column j, then denote λ by* ⟨i, j⟩;

 (3) *three beads of weight 1 on column(s) i, j and l, then denote λ by* ⟨i, j, l⟩.

*Remark.* 15 for all *p*-regular partitions of the above ⟨ ⟩-notation see Appendix A.1.

**Definition 3.2.** *let $\tilde{\lambda}$ be the partition corresponding to the Specht module belonging to $B_s$, it is clear that $\tilde{\lambda}$ can be represented on an abacus with 15 beads, and we make the following notation: if the abacus display of $\tilde{\lambda}$ with 15 beads has*

 (1) *one bead of weight 2 on column i, then denote $\tilde{\lambda}$ by* ⟨i⟩;

 (2) *two beads of weight 1 on column(s) i and j, then denote $\tilde{\lambda}$ by* ⟨i, j⟩.

*Remark.* The above ⟨ ⟩-notation corresponding to the partition of Specht modules belonging to $B_s$ is given in Appendix A.3-A.6.

**Definition 3.3.** *let $\tilde{\lambda}$ be the partition corresponding to the Specht module belonging to $B_1$, it is clear that $\tilde{\lambda}$ can be represented on an abacus with 15 beads, and we make the following notation: if the abacus display of $\tilde{\lambda}$ with 15 beads has one bead of weight 1 on column i, then denote $\tilde{\lambda}$ by* ⟨i⟩.

*Remark.* The above ⟨ ⟩-notation corresponding to the partition of Specht modules belonging to $B_s$ is given in Appendix A.2.

**Lemma 3.4.** *The simple modes of B are related to the tensor of the signature representation sgn as follows*:

$$D^{\langle 5 \rangle} \leftrightarrow D^{\langle 4,3,2 \rangle}; \quad D^{\langle 4 \rangle} \leftrightarrow D^{\langle 4,3,1 \rangle}; \quad D^{\langle 3 \rangle} \leftrightarrow D^{\langle 4,2,1 \rangle}; \quad D^{\langle 2 \rangle} \leftrightarrow D^{\langle 5,2,1 \rangle};$$

$D^{\langle 1 \rangle} \leftrightarrow D^{\langle 1,2 \rangle}$; $\quad D^{\langle 5,4 \rangle} \leftrightarrow D^{\langle 3,3,2 \rangle}$; $\quad D^{\langle 5,3 \rangle} \leftrightarrow D^{\langle 4,4,2 \rangle}$; $\quad D^{\langle 5,2 \rangle} \leftrightarrow D^{\langle 5,5,2 \rangle}$;

$D^{\langle 5,1 \rangle} \leftrightarrow D^{\langle 2,1 \rangle}$; $\quad D^{\langle 4,5 \rangle} \leftrightarrow D^{\langle 5,4,3 \rangle}$; $\quad D^{\langle 4,3 \rangle} \leftrightarrow D^{\langle 4,4,3 \rangle}$; $\quad D^{\langle 4,2 \rangle} \leftrightarrow D^{\langle 5,5,3 \rangle}$;

$D^{\langle 4,1 \rangle} \leftrightarrow D^{\langle 3,1 \rangle}$; $\quad D^{\langle 3,5 \rangle} \leftrightarrow D^{\langle 5,4,2 \rangle}$; $\quad D^{\langle 3,4 \rangle} \leftrightarrow D^{\langle 5,4,1 \rangle}$; $\quad D^{\langle 3,2 \rangle} \leftrightarrow D^{\langle 5,5,4 \rangle}$;

$D^{\langle 2,5 \rangle} \leftrightarrow D^{\langle 5,3,2 \rangle}$; $\quad D^{\langle 2,4 \rangle} \leftrightarrow D^{\langle 5,3,1 \rangle}$; $\quad D^{\langle 1,5 \rangle} \leftrightarrow D^{\langle 2,3 \rangle}$; $\quad D^{\langle 1,4 \rangle} \leftrightarrow D^{\langle 1,3 \rangle}$.

*Proof.* We have easy access by Lemma 2.4.

*Example.* Let $P=(15)$ be a 5-regular partition of 15, and solve for a 5-regular partition $Q$ of 15 such that $D^P \otimes \mathrm{sgn} = D^Q$, where sgn is the signature representation of $kS_{15}$. From Definition 2.2, we know that $P'=I(P)=(10)$, $P''=I(P')=(5)$, and $I(P'')= \emptyset$. So $\alpha=3$, $P_1 =(15)$, $P_2 =(10)$, $P_3 =(5)$; $a_1 =e(P_1)$, $a_2 =e(P_2)$, $a_3 =e(P_3)=5$; $r_1 = r_2 = r_3 =1$; $e_1 = e_2 = e_3 =0$; $s_1 = s_2 = s_3 =4$. It is also known from Lemma 2.3 that $(5,5,5)$ and $(4,4,4)$ uniquely determine a $p$-regular partition sequence $(Q_1,Q_2,Q_3)$ such that $Q_i$ has 5 rows and $e(Q_i)=5$. It is not difficult to find $Q_3 =(4,3)$, i.e., $Q =(4,3)$.

**Proposition 3.5** (Restriction). *Restriction of simple of B it to $S_{14}$ produces the following modules:*

$D^{\langle 5 \rangle} \downarrow_{S_{14}} = D^{(14)}$; $\qquad D^{\langle 4 \rangle} \downarrow_{S_{14}} = D^{(13,1)}$;

$D^{\langle 3 \rangle} \downarrow_{S_{14}} = D^{(12,1^2)}$; $\qquad D^{\langle 2 \rangle} \downarrow_{S_{14}} = D^{(11,1^3)}$;

$D^{\langle 1 \rangle} \downarrow_{S_{14}} = D^{(10,1^4)}$; $\qquad D^{\langle 5,4 \rangle} \downarrow_{S_{14}} = D^{(9,5)} + D^{(10,4)}$;

$D^{\langle 5,3 \rangle} \downarrow_{S_{14}} = D^{(9,4,1)} + D^{(10,3,1)}$; $\qquad D^{\langle 5,2 \rangle} \downarrow_{S_{14}} = D^{(9,3,1^2)} + D^{(10,2,1^2)}$;

$D^{\langle 5,1 \rangle} \downarrow_{S_{14}} = D^{(9,2,1^3)}$; $\qquad D^{\langle 4,5 \rangle} \downarrow_{S_{14}} = D^{(8,6)}$;

$D^{\langle 4,3 \rangle} \downarrow_{S_{14}} = D^{(8,4,2)} + D^{(9,3,2)}$; $\qquad D^{\langle 4,2 \rangle} \downarrow_{S_{14}} = D^{(8,3,2,1)} + D^{(9,2,2,1)}$;

$D^{\langle 4,1 \rangle} \downarrow_{S_{14}} = D^{(8,2^2,1^2)}$; $\qquad D^{\langle 3,5 \rangle} \downarrow_{S_{14}} = D^{(7,6,1)} + D^{(8,5,1)}$;

$D^{\langle 3,4 \rangle} \downarrow_{S_{14}} = D^{(7,5,2)}$; $\qquad D^{\langle 3,2 \rangle} \downarrow_{S_{14}} = D^{(7,3,2^2)} + D^{(8,2^3)}$;

$D^{\langle 3,1 \rangle} \downarrow_{S_{14}} = D^{(7,2^3,1)}$; $\qquad D^{\langle 2,5 \rangle} \downarrow_{S_{14}} = D^{(6^2,1^2)} + D^{(7,5,1^2)}$;

$D^{\langle 2,4 \rangle}\downarrow_{S_{14}} = D^{(6,5,2,1)} + D^{(7,4,2,1)}$;  $\quad D^{\langle 2,3 \rangle}\downarrow_{S_{14}} = D^{(6,4,2^2)}$;

$D^{\langle 2,1 \rangle}\downarrow_{S_{14}} = D^{(6,2^4)}$;  $\quad D^{\langle 1,5 \rangle}\downarrow_{S_{14}} = D^{(6,5,1^3)}$;

$D^{\langle 1,4 \rangle}\downarrow_{S_{14}} = D^{(5^2,2,1^2)} + D^{(6,4,2,1^2)}$;  $\quad D^{\langle 1,3 \rangle}\downarrow_{S_{14}} = D^{(5,4,2^2,1)} + D^{(6,3,2^2,1)}$;

$D^{\langle 1,2 \rangle}\downarrow_{S_{14}} = D^{(5,3,2^3)}$;  $\quad D^{\langle 5,4,3 \rangle}\downarrow_{S_{14}} = D^{(5^2,4)}$;

$D^{\langle 5,4,2 \rangle}\downarrow_{S_{14}} = D^{(5,4^2,1)} + D^{(5^2,3,1)}$;  $\quad D^{\langle 5,4,1 \rangle}\downarrow_{S_{14}} = D^{(5,4,3,1^2)}$;

$D^{\langle 5,5,4 \rangle}\downarrow_{S_{14}} = D^{(5,4,2,1^3)} + D^{(5^2,1^4)}$;  $\quad D^{\langle 5,3,2 \rangle}\downarrow_{S_{14}} = D^{(4^3,2)} + D^{(5,4,3,2)}$;

$D^{\langle 5,3,1 \rangle}\downarrow_{S_{14}} = D^{(4^2,3,2,1)} + D^{(5,3^2,2,1)}$;  $\quad D^{\langle 5,5,3 \rangle}\downarrow_{S_{14}} = D^{(4^2,2^2,1^2)} + D^{(5,3,2^2,1^2)}$;

$D^{\langle 5,2,1 \rangle}\downarrow_{S_{14}} = D^{(4,3^2,2^2)}$;  $\quad D^{\langle 5,5,2 \rangle}\downarrow_{S_{14}} = D^{(4,3,2^3,1)} + D^{(5,2^4,1)}$;

$D^{\langle 4,3,2 \rangle}\downarrow_{S_{14}} = D^{(4^2,3^2)}$;  $\quad D^{\langle 4,3,1 \rangle}\downarrow_{S_{14}} = D^{(4,3^3,1)}$;

$D^{\langle 4,4,3 \rangle}\downarrow_{S_{14}} = D^{(4,3,2^2,1^3)} + D^{(4^2,2,1^4)}$;  $\quad D^{\langle 4,2,1 \rangle}\downarrow_{S_{14}} = D^{(3^4,2)}$;

$D^{\langle 4,4,2 \rangle}\downarrow_{S_{14}} = D^{(3^2,2^3,1^2)} + D^{(4,2^4,1^2)}$;  $\quad D^{\langle 3,3,2 \rangle}\downarrow_{S_{14}} = D^{(3,2^4,1^3)} + D^{(3^2,2^2,1^4)}$.

*Proof.* We have easy access by Lemma 2.8.

*Example.* Consider

$$D^{\langle 5,3 \rangle}\downarrow_{S_{14}}, \text{where } \langle 5,3 \rangle = (10,4,1).$$

By Lemma 2.7, 1 is normal; when $i=2$, since $\beta(1,2)=8\not\equiv 0 \pmod 5$, 2 is normal. When $i=3$, since $\beta(i-1,i)=\beta(2,3)=5\equiv 0\pmod 5$, 3 is not normal. Thus, only 1 and 2 are normal. Since $\lambda(1)=(9,4,1)$, $\lambda(2)=(10,3,1)$, by Lemma 2.8, we have

$$D^{\langle 5,3 \rangle}\downarrow_{S_{14}} = D^{(9,4,1)} + D^{(10,3,1)}.$$

**Lemma 3.6.** *If $\{\tilde{\lambda}_1, \tilde{\lambda}_2, \ldots, \tilde{\lambda}_a\}$ is the complete set of partitions corresponding to Specht modules all belonging to $B_s$, and we denote these partitions by the $\langle \ \rangle$- notation of the abacus display with 15* beads, then $\{\tilde{\lambda}_1,\tilde{\lambda}_2,\ldots,\tilde{\lambda}_a\}$ consists of $\langle t,u\rangle$, $\langle u\rangle$, $\langle s-1,t\rangle$, $\langle s,t\rangle$, $\langle s-1,s-1\rangle$, $\langle s-1\rangle$, $\langle s,s\rangle$, $\langle s\rangle$, $\langle s-1,s\rangle$, where $t,u\neq s-1,s$, and there is:

$$S^{\langle t,u\rangle}\uparrow^B \sim S^{\langle s,t,u\rangle} + S^{\langle s-1,t,u\rangle}; \qquad S^{\langle u\rangle}\uparrow^B \sim S^{\langle u,s\rangle} + S^{\langle u,s-1\rangle};$$

$$S^{\langle s-1,t\rangle}\uparrow^B \sim S^{\langle s,t\rangle} + S^{\langle s-1,t\rangle}; \qquad S^{\langle s,t\rangle}\uparrow^B \sim S^{\langle s,s,t\rangle} + S^{\langle s-1,s-1,t\rangle}.$$

*Proof.* We have easy access by Lemma 2.10.

**Proposition 3.7**(Induction 1). *Induction of simple modules on B in $B_5$ has the following Loewy modular structure:*

$$D^{\langle 4\rangle}\uparrow^B = \begin{matrix}D^{\langle 5\rangle}\\D^{\langle 4\rangle}\\D^{\langle 5\rangle}\end{matrix}; \qquad D^{\langle 4,4\rangle}\uparrow^B = \begin{matrix}D^{\langle 5,4\rangle}\\D^{\langle 4,5\rangle}D^{\langle 4\rangle}\\D^{\langle 5,4\rangle}\end{matrix}; \qquad D^{\langle 4,3\rangle}\uparrow^B = \begin{matrix}D^{\langle 5,3\rangle}\\D^{\langle 4,3\rangle}D^{\langle 3\rangle}\\D^{\langle 5,3\rangle}\end{matrix};$$

$$D^{\langle 4,2\rangle}\uparrow^B = \begin{matrix}D^{\langle 5,2\rangle}\\D^{\langle 4,2\rangle}D^{\langle 2\rangle}\\D^{\langle 5,2\rangle}\end{matrix}; \qquad D^{\langle 4,1\rangle}\uparrow^B = \begin{matrix}D^{\langle 5,1\rangle}\\D^{\langle 4,1\rangle}D^{\langle 1\rangle}\\D^{\langle 5,1\rangle}\end{matrix}; \qquad D^{\langle 3\rangle}\uparrow^B = \begin{matrix}D^{\langle 3,5\rangle}\\D^{\langle 3,4\rangle}D^{\langle 4,5\rangle}\\D^{\langle 3,5\rangle}\end{matrix};$$

$$D^{\langle 2\rangle}\uparrow^B = \begin{matrix}D^{\langle 2,5\rangle}\\D^{\langle 2,4\rangle}\\D^{\langle 2,5\rangle}\end{matrix}; \qquad D^{\langle 1\rangle}\uparrow^B = \begin{matrix}D^{\langle 1,5\rangle}\\D^{\langle 1,4\rangle}\\D^{\langle 1,5\rangle}\end{matrix}; \qquad D^{\langle 5\rangle}\uparrow^B = \begin{matrix}D^{\langle 5,5,4\rangle}\\D^{\langle 4,1\rangle}\\D^{\langle 5,5,4\rangle}\end{matrix};$$

$$D^{\langle 3,2\rangle}\uparrow^B = \begin{matrix}D^{\langle 5,3,2\rangle}\\D^{\langle 2,3\rangle}D^{\langle 4,3,2\rangle}\\D^{\langle 5,3,2\rangle}\end{matrix}; \qquad D^{\langle 3,1\rangle}\uparrow^B = \begin{matrix}D^{\langle 5,3,1\rangle}\\D^{\langle 4,3,1\rangle}D^{\langle 1,3\rangle}\\D^{\langle 5,3,1\rangle}\end{matrix}; \qquad D^{\langle 5,3\rangle}\uparrow^B = \begin{matrix}D^{\langle 5,5,3\rangle}\\D^{\langle 4,4,3\rangle}D^{\langle 3,1\rangle}\\D^{\langle 5,5,3\rangle}\end{matrix};$$

$$D^{\langle 2,1\rangle}\uparrow^B = \begin{matrix}D^{\langle 5,2,1\rangle}\\D^{\langle 4,2,1\rangle}D^{\langle 4,3,1\rangle}\\D^{\langle 5,2,1\rangle}\end{matrix}; \qquad D^{\langle 5,2\rangle}\uparrow^B = \begin{matrix}D^{\langle 5,5,2\rangle}\\D^{\langle 4,4,2\rangle}D^{\langle 2,1\rangle}\\D^{\langle 5,5,2\rangle}\end{matrix}.$$

*Proof.* We consider the Loewy structure of the induction of simple modules in $B_5$ onto $B$ in order of lexicographic order from largest to smallest. Consider first $D^{\langle 4\rangle}\uparrow^B$, which clearly has $D^{\langle 4\rangle} = D^{(1^4)} = S^{(1^4)}$. By Lemma 3.6 we have

$$S^{(1^4)}\uparrow^B_{B_5} \sim S^{(1^5)} + S^{(1^4,1)}.$$

And by the decomposition matrix of $S_{15}(p=5)$ we have

$$S^{(1^5)} = S^{\langle 5\rangle} \sim D^{\langle 5\rangle}, \quad S^{(1^4,1)} = S^{\langle 4\rangle} \sim D^{\langle 5\rangle} + D^{\langle 4\rangle}.$$

so

$$D^{\langle 4\rangle}\uparrow^B = 2D^{\langle 5\rangle} + D^{\langle 4\rangle}.$$

By duality we have:

$$D^{\langle 4\rangle}\uparrow^B = \begin{matrix} D^{\langle 5\rangle} \\ D^{\langle 4\rangle} \\ D^{\langle 5\rangle} \end{matrix}.$$

Now for any simple module $D^{\tilde{\lambda}}$ in $B_5$, analyze the composition factor of $D^{\tilde{\lambda}}$. Assume that we already know the composition factor of $D^{\tilde{\mu}}$ for all partitions $\tilde{\mu}$ that satisfy $\tilde{\mu} > \tilde{\lambda}$ and whose corresponding Specht module belongs to $B_5$. First, the composition factor of $S^{\tilde{\lambda}}\uparrow^B$ can be obtained from Lemma 3.6 as well as the decomposition matrix of $S_{15}$, and then the composition factor of $S^{\tilde{\lambda}}$ can be known from the decomposition matrix of $S_{14}$. Taking any one of the composition factors of $S^{\tilde{\lambda}}$, $D^\rho$, and $\rho \neq \tilde{\lambda}$, it is not difficult to find that $\rho > \tilde{\lambda}$, and thus the composition factor of $D^\rho\uparrow^B$ can be known. Finally, subtracting the composition factor of $S^{\tilde{\lambda}}\uparrow^B$ from the composition factor of $D^\rho\uparrow^B (\rho \neq \tilde{\lambda})$, what we get is the composition factor of $D^{\tilde{\lambda}}\uparrow^B$, and then using the duality that is the Loewy structure of $D^{\tilde{\lambda}}\uparrow^B$.

**Lemma 3.8.** *The simple modules of B2 and B5 are related to the tensor of the signature representation sgn as follows:*

$D^{(1^4)} \leftrightarrow D^{(4^2,3^2)}$; $\quad D^{(9,5)} \leftrightarrow D^{(3,2^4,1^3)}$; $\quad D^{(9,4,1)} \leftrightarrow D^{(4,2^4,1^2)}$; $\quad D^{(9,3,1^2)} \leftrightarrow D^{(5,2^4,1)}$;

$D^{(9,2,1^3)} \leftrightarrow D^{(6,2^4)}$; $\quad D^{(8,5,1)} \leftrightarrow D^{(5^2,3,1)}$; $\quad D^{(7,5,1^2)} \leftrightarrow D^{(5,4,3,2)}$; $\quad D^{(6,5,1^3)} \leftrightarrow D^{(6,4,2^2)}$;

$D^{(5^2,1^4)} \leftrightarrow D^{(8,2^3)}$; $\quad D^{(4^3,2)} \leftrightarrow D^{(6^2,1^2)}$; $\quad D^{(4^2,3,2,1)} \leftrightarrow D^{(6,5,2,1)}$; $\quad D^{(4^2,2^2,1^2)} \leftrightarrow D^{(9,2^2,1)}$;

$D^{(4,3^2,2^2)} \leftrightarrow D^{(11,1^3)}$; $\quad D^{(4,3,2^3,1)} \leftrightarrow D^{(10,2,1^2)}$.

*Proof.* We have easy access by Lemma 2.4.

*Remark.* In fact, the tensor of the signature representation sgn is related to the simple modules in $B_4$ and $B_3$ in a similar way. That is, given any partition $\tilde{\lambda}$ that corresponds to a Specht module in $B_4$, there is a uniquely determined partition $\tilde{\mu}$ such that $D^{\tilde{\lambda}} \otimes \text{sgn} = D^{\tilde{\mu}}$ and the Specht module associated with $\tilde{\mu}$ is in $B_3$. Moreover, it is not difficult to see that for the two partitions that have the above for relation, their corresponding simple modules induced on B have a similar correspondence, i.e., $D^{\tilde{\mu}}\uparrow^B = D^{\tilde{\lambda}}\uparrow^B \otimes \text{sgn}$.

**Lemma 3.9.** *If $\{\tilde{\lambda}_1,\tilde{\lambda}_2,\ldots,\tilde{\lambda}_a\}$ is the complete set of partitions corresponding to Specht modules all belonging to $B_1$, and we denote these partitions by the $\langle\,\rangle$- notation of the abacus display with 15 beads, then $\{\tilde{\lambda}_1,\tilde{\lambda}_2,\ldots,\tilde{\lambda}_a\}$ consists of $\langle 1,1\rangle$, $\langle 1,i\rangle$, $\langle 1\rangle$, where $2\leq i\leq p-1$. So following the form*

$$S^{\tilde{\lambda}}\uparrow_{B_1}^{B} \sim S^{\lambda} + S^{\mu} + S^{\upsilon}, \ \lambda > \mu > \upsilon,$$

*we have:*

$$S^{\langle 1 \rangle}\uparrow^{B} \sim S^{\langle 1,1 \rangle} + S^{\langle 1,1,p \rangle} + S^{\langle p,p,p \rangle};$$

$$S^{\langle i \rangle}\uparrow^{B} \sim S^{\langle 1,i \rangle} + S^{\langle 1,i,p \rangle} + S^{\langle i,p,p \rangle};$$

$$S^{\langle p \rangle}\uparrow^{B} \sim S^{\langle 1 \rangle} + S^{\langle p,1 \rangle} + S^{\langle p,p \rangle}.$$

*Proof.* We have easy access by Lemma 2.14.

**Proposition 3.10**(Induction 2). *Induction of simple modules on B in $B_1$ has the following Loewy modular structure:*

$$D^{\langle 5 \rangle}\uparrow^{B} = \begin{matrix} D^{\langle 1 \rangle} \\ D^{\langle 5,1 \rangle} D^{\langle 2 \rangle} \\ D^{\langle 5,2 \rangle} D^{\langle 1 \rangle} \\ D^{\langle 5,1 \rangle} D^{\langle 2 \rangle} \\ D^{\langle 1 \rangle} \end{matrix} ; \qquad D^{\langle 4 \rangle}\uparrow^{B} = \begin{matrix} D^{\langle 1,4 \rangle} \\ D^{\langle 5,4 \rangle} D^{\langle 2,4 \rangle} \\ D^{\langle 5,5,4 \rangle} D^{\langle 4,2 \rangle} D^{\langle 1,4 \rangle} \\ D^{\langle 5,4 \rangle} D^{\langle 2,4 \rangle} \\ D^{\langle 1,4 \rangle} \end{matrix} ;$$

$$D^{\langle 3 \rangle}\uparrow^{B} = \begin{matrix} D^{\langle 1,3 \rangle} \\ D^{\langle 5,3,1 \rangle} D^{\langle 3,4 \rangle} \\ D^{\langle 5,5,3 \rangle} D^{\langle 1,3 \rangle} D^{\langle 3,2 \rangle} \\ D^{\langle 5,3,1 \rangle} D^{\langle 3,4 \rangle} \\ D^{\langle 1,3 \rangle} \end{matrix} ; \qquad D^{\langle 2 \rangle}\uparrow^{B} = \begin{matrix} D^{\langle 1,2 \rangle} \\ D^{\langle 5,2,1 \rangle} D^{\langle 2,1 \rangle} \\ D^{\langle 5,5,2 \rangle} D^{\langle 1,2 \rangle} \\ D^{\langle 5,2,1 \rangle} D^{\langle 2,1 \rangle} \\ D^{\langle 1,2 \rangle} \end{matrix} .$$

*Proof.* The proof here is similar to that of Proposition 3.7.

**Proposition 3.11.** *[1, Proposition 2.10] Let $D^{\lambda}$ and $D^{\mu}$ be simple modules in B. Then $Ext^1(D^{\lambda},D^{\mu}) \neq 0$ holds if and only if one of the following two conditions holds:*

*(1) $D^{\mu}$ lies in the second layer of the Loewy structure of $D^{\tilde{\lambda}}\uparrow^{B}$, where $D^{\tilde{\lambda}}$ is a direct summand of $D^{\mu}\downarrow_{S_{14}}$;*

*(2) There exist $D^{\tilde{\lambda}}$ and $D^{\tilde{\mu}}$ direct summands $D^{\lambda}\downarrow_{S_{14}}$ and $D^{\mu}\downarrow_{S_{14}}$ respectively such that $Ext^1(D^{\lambda},D^{\mu}) \neq 0$.*

**Theorem 3.12.** *Let $\lambda$ and $\mu$ be 5-regular partitions of 15 corresponding to Specht modules both belonging to B and $\lambda > \mu$. Then the following is a list of complete sets of $\lambda$ and $\mu$ such that $Ext^1(D^{\lambda},D^{\mu}) \neq 0$.*

| $\lambda$ | $\mu$ |
|:---:|:---:|
| ⟨5⟩ | ⟨4⟩,⟨5,3⟩,⟨2,5⟩ |
| ⟨4⟩ | ⟨3⟩,⟨5,4⟩,⟨4,3⟩,⟨2,4⟩ |
| ⟨3⟩ | ⟨2⟩,⟨5,3⟩,⟨3,4⟩ |
| ⟨2⟩ | ⟨1⟩,⟨5,2⟩,⟨2,4⟩ |
| ⟨1⟩ | ⟨5,1⟩,⟨1,4⟩ |
| ⟨5,4⟩ | ⟨5,3⟩,⟨4,5⟩ |
| ⟨5,3⟩ | ⟨5,2⟩,⟨4,3⟩,⟨3,5⟩ |
| ⟨5,2⟩ | ⟨5,1⟩,⟨4,2⟩,⟨2,5⟩ |
| ⟨5,1⟩ | ⟨4,1⟩,⟨1,5⟩ |
| ⟨4,5⟩ | ⟨4,3⟩,⟨3,5⟩,⟨2,4⟩,⟨5,4,1⟩ |
| ⟨4,3⟩ | ⟨4,2⟩,⟨3,4⟩ |
| ⟨4,2⟩ | ⟨4,1⟩,⟨3,2⟩,⟨2,4⟩ |
| ⟨4,1⟩ | ⟨3,1⟩,⟨1,4⟩,⟨5,5,4⟩ |
| ⟨3,5⟩ | ⟨3,4⟩,⟨2,5⟩,⟨5,3,1⟩ |
| ⟨3,4⟩ | ⟨3,2⟩,⟨2,4⟩,⟨1,3⟩,⟨5,4,3⟩ |
| ⟨3,2⟩ | ⟨3,1⟩,⟨2,3⟩ |
| ⟨3,1⟩ | ⟨2,1⟩,⟨1,3⟩,⟨5,5,3⟩ |
| ⟨2,5⟩ | ⟨2,4⟩,⟨1,5⟩,⟨5,3,2⟩ |
| ⟨2,4⟩ | ⟨2,3⟩,⟨1,4⟩,⟨5,4,2⟩ |
| ⟨2,3⟩ | ⟨2,1⟩,⟨1,3⟩,⟨5,3,2⟩ |
| ⟨2,1⟩ | ⟨1,2⟩,⟨5,5,2⟩ |
| ⟨1,5⟩ | ⟨1,4⟩,⟨5,5,4⟩,⟨5,3,1⟩ |
| ⟨1,4⟩ | ⟨1,3⟩,⟨5,4,1⟩ |
| ⟨1,3⟩ | ⟨1,2⟩,⟨5,3,1⟩ |
| ⟨1,2⟩ | ⟨5,2,1⟩ |
| ⟨5,5,4⟩ | ⟨5,5,3⟩ |
| ⟨5,5,3⟩ | ⟨5,5,2⟩,⟨4,4,3⟩ |
| ⟨5,5,2⟩ | 0 |
| ⟨5,4,3⟩ | ⟨5,4,2⟩,⟨5,3,1⟩,⟨4,4,3⟩,⟨3,3,2⟩ |
| ⟨5,4,2⟩ | ⟨5,4,1⟩,⟨5,3,2⟩,⟨4,4,2⟩ |
| ⟨5,4,1⟩ | ⟨5,5,4⟩,⟨5,3,1⟩,⟨4,4,3⟩,⟨4,1,1⟩ |
| ⟨5,3,2⟩ | ⟨5,3,1⟩,⟨5,5,2⟩,⟨4,3,2⟩ |
| ⟨5,3,1⟩ | ⟨5,5,3⟩,⟨5,2,1⟩,⟨4,3,1⟩ |
| ⟨5,2,1⟩ | ⟨5,5,2⟩,⟨4,2,1⟩ |
| ⟨4,4,3⟩ | ⟨4,4,2⟩ |
| ⟨4,4,2⟩ | 0 |
| ⟨4,3,2⟩ | ⟨4,3,1⟩,⟨4,4,2⟩ |
| ⟨4,3,1⟩ | ⟨4,4,3⟩,⟨4,2,1⟩,⟨3,3,2⟩ |
| ⟨4,2,1⟩ | ⟨4,4,2⟩ |
| ⟨3,3,2⟩ | 0 |

*Proof.* We have this list by Proposition 3.5, Proposition 3.7, Lemma 3.8, Proposition 3.10, Lemma 2.18, Proposition 3.11, and the structure of the Ext-quiver of $B_i (i \in \{1,2,3,4,5\})$.

Thus we can obtain the Ext-quiver structure of $B$.

**Appendix A**

The ⟨ ⟩-notation for partitions represented by 15 beads on abacus.

**A.1.** The principal block $B(kS_{15})$

(15)= ⟨5⟩,                (14,1)= ⟨4⟩,              (13,1²)= ⟨3⟩,              (12,1³)= ⟨2⟩,

(11,1⁴)= ⟨1⟩,             (10,5)= ⟨5,4⟩,            (10,4,1)= ⟨5,3⟩,           (10,3,1²)= ⟨5,2⟩,

(10,2,1³)= ⟨5,1⟩,         (9,6)= ⟨4,5⟩,             (9,4,2)= ⟨4,3⟩,            (9,3,2,1)= ⟨4,2⟩,

(9,2²,1²)= ⟨4,1⟩,         (8,6,1)= ⟨3,5⟩,           (8,5,2)= ⟨3,4⟩,            (8,3,2²)= ⟨3,2⟩,

(8,2³, 1)= ⟨3,1⟩,         (7,6,1²)= ⟨2,5⟩,          (7,5,2,1)= ⟨2,4⟩,          (7,4,2²)= ⟨2,3⟩,

(8,2⁴)= ⟨2,1⟩,            (6²,1³)= ⟨1,5⟩,           (6,5,2,1²)= ⟨1,4⟩,         (6,4,2², 1)= ⟨1,3⟩,

(6,3,2³)= ⟨1,2⟩,          (5²,2,1³)= ⟨5,5,4⟩,       (5,4,2², 1²)= ⟨5,5,3⟩,     (5,3,2³, 1)= ⟨5,5,2⟩,

(5³)= ⟨5,4,3⟩,            (5², 4,1)= ⟨5,4,2⟩,       (5², 3,1²)= ⟨5,4,1⟩,       (5,4²·2)= ⟨5,3,2⟩,

(5,4,3,2,1)= ⟨5,3,1⟩,     (5,3²,2²)= ⟨5,2,1⟩,       (4²,2²,1³)= ⟨4,4,3⟩,       (4,3,2³, 1²)= ⟨4,4,2⟩,

(4³, 3)= ⟨4,3,2⟩,         (4², 3², 1)= ⟨4,3,1⟩,     (4, 3³, 2)= ⟨4,2,1⟩,       (3², 2³, 1³)= ⟨3,3,2⟩.

**A.2.** The principal block $B_1(kS_{14})$

(10,1⁴)= ⟨5⟩,     (5²,2,1²)= ⟨4⟩,     (5,4,2²,1)= ⟨3⟩,     (5,3,2³)= ⟨2⟩.

**A.3.** The principal block $B_2(kS_{14})$

(11,1³)= ⟨1⟩,         (10,2,1²)= ⟨5⟩,       (9,2²,1)= ⟨4⟩,        (8,2³)= ⟨3⟩,

(6²,1²)= ⟨5,1⟩,       (6,5,2,1)= ⟨4,1⟩,     (6,4,2²)= ⟨3,1⟩,      (6,2⁴)= ⟨1,1⟩,

(5²,3,1)= ⟨5,4⟩,      (5,4,3,2)= ⟨5,3⟩,     (5,2⁴, 1)= ⟨5,5⟩,     (4²,3²)= ⟨4,3⟩,

(4,2⁴, 1²)= ⟨4,4⟩,    (3,2⁴, 1³)= ⟨3,3⟩.

**A.4.** The principal block $B_3(kS_{14})$

(12,1²)= ⟨2⟩,         (10,3,1)= ⟨5⟩,        (9,3,2)= ⟨4⟩,         (7,6,1)= ⟨5,2⟩,

(7,5,2)= ⟨4,2⟩,       (7,3,2²)= ⟨2,2⟩,      (7,2³, 1)= ⟨2,1⟩,     (6,3,2², 1)= ⟨1⟩,

$(5^2,4)= \langle 5,4 \rangle$,  $(5,3^2,2,1)= \langle 5,1 \rangle$,  $(5,3,2^2,1^2)= \langle 5,5 \rangle$,  $(4,3^3,1)= \langle 4,1 \rangle$,

$(4,3,2^2,1^3)= \langle 4,4 \rangle$,  $(3^2,2^2,1^4)= \langle 3 \rangle$.

**A.5.** The principal block $B_4(kS_{14})$

$(13,1)= \langle 3 \rangle$,  $(10,4)= \langle 5 \rangle$,  $(8,6)= \langle 5,3 \rangle$,  $(8,4,2)= \langle 3,3 \rangle$,

$(8,3,2,1)= \langle 3,2 \rangle$,  $(8,2^2,1^2)= \langle 3,1 \rangle$,  $(7,4,2,1)= \langle 2 \rangle$,  $(6,4,2,1^2)= \langle 1 \rangle$,

$(5,4^2,1)= \langle 5,2 \rangle$,  $(5,4,3,1^2)= \langle 5,1 \rangle$,  $(5,4,2,1^3)= \langle 5,5 \rangle$,  $(4^2,2,1^4)= \langle 4 \rangle$,

$(3^4,2)= \langle 2,1 \rangle$,  $(3^2,2^3,1^2)= \langle 4,2 \rangle$.

**A.6.** The principal block $B_5(kS_{14})$

$(14)= \langle 4 \rangle$,  $(9,5)= \langle 4,4 \rangle$,  $(9,4,1)= \langle 4,3 \rangle$,  $(9,3,1^2)= \langle 4,2 \rangle$,

$(9,2,1^3)= \langle 4,1 \rangle$,  $(8,5,1)= \langle 3 \rangle$,  $(7,5,1^2)= \langle 2 \rangle$,  $(6,5,1^3)= \langle 1 \rangle$,

$(5^2,1^4)= \langle 5 \rangle$,  $(4^3,2)= \langle 3,2 \rangle$,  $(4^2,3,2,1)= \langle 3,1 \rangle$,  $(4^2,2^2,1^2)= \langle 5,3 \rangle$,

$(4,3^2,2^2)= \langle 2,1 \rangle$,  $(4,3,2^3,1)= \langle 5,2 \rangle$.

## Appendix B

The Ext-quiver of five block of $kS_{14}$.

**B.1.** The principal block $B_1$

| $\lambda$ | $\mu$ |
| --- | --- |
| $(10,1^4)$ | $(5^2,2,1^2)$ |
| $(5^2,2,1^2)$ | $(5,4,2^2,1)$ |
| $(5,4,2^2,1)$ | $(5,3,2^3)$ |

**B.2.** The principal block $B_2$

| $\lambda$ | $\mu$ |
| --- | --- |
| $(11,1^3)$ | $(10,2,1^2), (6,5,2,1)$ |
| $(10,2,1^2)$ | $(9,2^2,1), (6^2,1^2)$ |
| $(9,2^2,1)$ | $(8,2^3), (6,5,2,1)$ |
| $(8,2^3)$ | $(6,4,2^2)$ |
| $(6^2,1^2)$ | $(6,5,2,1), (5,4,3,2)$ |
| $(6,5,2,1)$ | $(6,4,2^2), (5^2,3,1)$ |
| $(6,4,2^2)$ | $(6,2^4), (5,4,3,2)$ |
| $(6,2^4)$ | $(5,2^4,1)$ |
| $(5^2,3,1)$ | $(5,4,3,2), (4,2^4,1^2)$ |
| $(5,4,3,2)$ | $(5,2^4,1), (4^2,3^2)$ |
| $(5,2^4,1)$ | $(4,2^4,1^2)$ |
| $(4^2,3^2)$ | $(4,2^4,1^2)$ |
| $(4,2^4,1^2)$ | $(3,2^4,1^3)$ |

**B.3.** The principal block $B_3$

| $\lambda$ | $\mu$ |
|---|---|
| $(12,1^2)$ | $(10,3,1), (7,5,2)$ |
| $(10,3,1)$ | $(9,3,2), (7,6,1)$ |
| $(9,3,2)$ | $(7,5,2)$ |
| $(7,6,1)$ | $(7,5,2),(5,3^2,2,1)$ |
| $(7,5,2)$ | $(7,3,2^2), (6,3,2^2,1)$ |
| $(7,3,2^2)$ | $(7,2^3,1)$ |
| $(7,2^3,1)$ | $(6,3,2^2,1),(5,3,2^2,1^2)$ |
| $(6,3,2^2,1)$ | $(5,3^2,2,1)$ |
| $(5^2,4)$ | $(5,3^2,2,1),(4,3,2^2,1^3),(3^2,2^2,1^4)$ |
| $(5,3^2,2,1)$ | $(5,3,2^2,1^2),(4,3^3,1)$ |
| $(5,3,2^2,1^2)$ | $(4,3,2^2,1^3)$ |
| $(4,3^3,1)$ | $(4,3,2^2,1^3),(3^2,2^2,1^4)$ |

**B.4.** The principal block $B_4$

| $\lambda$ | $\mu$ |
|---|---|
| $(13,1)$ | $(10,4),(8,4,2),(7,4,2,1)$ |
| $(10,4)$ | $(8,6)$ |
| $(8,6)$ | $(8,4,2),(7,4,2,1),(5,4,3,1^2)$ |
| $(8,4,2)$ | $(8,3,2,1)$ |
| $(8,3,2,1)$ | $(8,2^2,1^2),(7,4,2,1)$ |
| $(8,2^2,1^2)$ | $(6,4,2,1^2),(5,4,2,1^3)$ |
| $(7,4,2,1)$ | $(6,4,2,1^2),(5,4^2,1)$ |
| $(6,4,2,1^2)$ | $(5,4,3,1^2)$ |
| $(5,4^2,1)$ | $(5,4,3,1^2),(3^2,2^3,1^2)$ |
| $(5,4,3,1^2)$ | $(5,4,2,1^3),(4^2,2,1^4),(3^4,2)$ |
| $(4^2,2,1^4)$ | $(3^2,2^3,1^2)$ |
| $(3^4,2)$ | $(3^2,2^3,1^2)$ |

**B.5.** The principal block $B_5$

| $\lambda$ | $\mu$ |
|---|---|
| $(14)$ | $(9,4,1),(7,5,1^2)$ |
| $(9,5)$ | $(9,4,1)$ |
| $(9,4,1)$ | $(9,3,1^2),(8,5,1)$ |
| $(9,3,1^2)$ | $(9,2,1^3),(7,5,1^2)$ |
| $(9,2,1^3)$ | $(6,5,1^3)$ |
| $(8,5,1)$ | $(7,5,1^2), (4^2,3,2,1)$ |
| $(7,5,1^2)$ | $(6,5,1^3), (4^3,2)$ |
| $(6,5,1^3)$ | $(5^2,1^4), (4^2,3,2,1)$ |
| $(5^2,1^4)$ | $(4^2,2^2,1^2)$ |
| $(4^3,2)$ | $(4^2,3,2,1), (4,3,2^3,1)$ |
| $(4^2,3,2,1)$ | $(4^2,2^2,1^2), (4,3^2,2^2)$ |
| $(4^2,2^2,1^2)$ | $(4,3,2^3,1)$ |
| $(4,3^2,2^2)$ | $(4,3,2^3,1)$ |